\documentclass[11pt,a4paper]{article}
\usepackage{a4wide}

\usepackage[T1]{fontenc}
\usepackage{lmodern}
\usepackage{amsthm,amsmath,amssymb}
\usepackage{hyperref}
\usepackage{microtype}
\usepackage{hyperref}
\usepackage{graphicx}
\usepackage{enumitem}
\usepackage{fdsymbol}
\usepackage{bbding}
\usepackage{fontawesome}
\usepackage{changepage}

\newcommand\cB{{\mathcal B}}
\newcommand\cC{{\mathcal C}}

\newcommand\cF{{\mathcal F}}
\newcommand\cG{{\mathcal G}}

\newcommand\cN{{\mathcal N}}

\makeatletter
\newtheorem*{rep@thm}{\rep@title}
\newcommand{\newreptheorem}[2]{%
\newenvironment{rep#1}[1]{%
 \def\rep@title{#2 \ref{##1}}%
 \begin{rep@thm}}%
 {\end{rep@thm}}}
\makeatother

\theoremstyle{plain}
\newtheorem{thm}{Theorem}[section]
\newtheorem{theorem}{Theorem}

\newreptheorem{thm}{Theorem}
\newtheorem{lemma}[thm]{Lemma}

\newtheorem{proposition}[thm]{Proposition}
\newtheorem{observation}[thm]{Observation}

\theoremstyle{definition}

\newtheorem{claim}[thm]{Claim}

\title{The Constructor-Blocker Game}
\author{Bal\'azs Patk\'os\footnote{Alfr\'ed R\'enyi Institute of Mathematics, Budapest, Hungary. Partially supported by NKFIH grants SNN 129364 and FK 132060.} \and Milo\v{s} Stojakovi\'{c}\footnote{Department of Mathematics and Informatics, Faculty of Sciences, University of Novi Sad, Serbia. Partly supported by Ministry of Education, Science and Technological Development of the Republic of Serbia (Grant No.~451-03-68/2022-14/200125). Partly supported by Provincial Secretariat for Higher Education and Scientific Research, Province of Vojvodina (Grant No.~142-451-2686/2021). {\tt milos.stojakovic@dmi.uns.ac.rs}} \and M\'at\'e Vizer \footnote{Alfr\'ed R\'enyi Institute of Mathematics, Budapest, Hungary. Partially supported by NKFIH grants SNN 129364 and FK 132060  by the  J\'anos Bolyai Research Fellowship  and by the New National Excellence Program under the grant number \'UNKP-21-5-BME-361. {\tt vizermate@gmail.com}}
}

\date{}

\begin{document}

\maketitle

\begin{abstract}
    We study the following game version of generalized graph Tur\'an problems. For two fixed graphs $F$ and $H$, two players, Constructor and Blocker, alternately claim  unclaimed edges of the complete graph $K_n$. Constructor can only claim edges so that he never claims all edges of any copy of $F$, i.e.~his graph must remain $F$-free, while Blocker can claim unclaimed edges without restrictions. The game ends when Constructor cannot claim further edges or when all edges have been claimed. The score of the game is the number of copies of $H$ with all edges claimed by Constructor. Constructor's aim is to maximize the score, while Blocker tries to keep the score as low as possible.
We denote by $g(n,H,F)$ the score of the game when both players play optimally and Constructor starts the game.

In this paper, we obtain the exact value of $g(n,H,F)$ when both $F$ and $H$ are stars and when $F=P_4$, $H=P_3$. We determine the asymptotics of $g(n,H,F)$ when $F$ is a star and $H$ is a tree and when $F=P_5$, $H=K_3$, and we derive upper and lower bounds on $g(n,P_4,P_5)$.
\end{abstract}

\section{Introduction}

The Tur\'an problem for a set $\cF$ of graphs asks the following: What is the maximum number $ex(n,\cF)$ of edges that a graph on $n$ vertices can have without containing any $F \in \cF$ as a subgraph? When $\cF$ contains a single graph $F$, we simply  write $ex(n,F)$. This function has been intensively studied, starting with Mantel \cite{M1907} and Tur\'an \cite{T1941} who determined $ex(n,K_r)$ where $K_r$ denotes the complete graph on $r$ vertices with $r \ge 3$. See \cite{FS2013,S1997} for surveys on this topic.

There is a generalization of the Tur\'an problem, when we count the maximum number of copies of a certain graph $H$ in a graph $G$ on $n$ vertices, provided that $G$ does not contain any $F \in \cF$ as a subgraph, for a set $\cF$ of graphs. To be more precise, let us introduce some notation: for two graphs $H$ and $G$, let $\cN(H,G)$ denote the number of copies of $H$ in $G$. We say that a graph $G$ is \textit{$\cF$-free}, if $G$ does not contain any $F \in \cF$ as a subgraph. Given a graph $H$ and a set $\cF$ of graphs, let $$ex(n,H,\cF)=\max_{G} \{\cN(H,G): \text{$G$ is an $\cF$-free graph on $n$ vertices} \}.$$
If $\cF=\{F\}$, we simply denote it by $ex(n,H,F)$. Note that the classical Tur\'an problem can be stated using this notion as $ex(n,F)=ex(n,K_2,F)$. This problem was initiated by Zykov~\cite{Z1949}, who determined $ex(n,K_s,K_t)$ exactly, however the systematic study of the function $ex(n,H,F)$ started just recently in~\cite{ALS2016} by Alon and Shikhelman. This topic became extensively investigated in the recent years, see e.g.~\cite{GySTZ19,GMV19,GGyMV20} and the references therein.

\medskip 

Our goal is to introduce a game analogue of the parameter $ex(n,H,F)$ and provide some results. For two fixed graphs $F$ and $H$, two players, Constructor and Blocker  alternately claim an unclaimed edge of the complete graph $K_n$. Constructor can only claim edges so that he never claims all edges of any copy of $F$, i.e.~his graph must remain $F$-free, while Blocker can claim unclaimed edges without restrictions. The game ends when Constructor cannot claim further edges or when all edges have been claimed. The score of the game is the number of copies of $H$ with all edges claimed by Constructor. Constructor's aim is to maximize the score, while Blocker tries to keep the score as low as possible. 
We denote by $g(n,H,F)$ the score of the game when both players play optimally and Constructor starts the game. Let us note that for all the games that we study in this paper the identity of the starting player does not have a big impact -- it turns out that all the results still hold if Blocker starts the game. 

Our Constructor-Blocker games borrow some aspects of two well-studied classes of combinatorial games on graphs -- the Maker-Breaker positional games, and the saturation games. It turns out that these two settings fit together well to give us solid ground for a game version of the above mentioned generalization of the Tur\'an problem.

In a Maker-Breaker game the players alternately claim unclaimed elements of the board $X$, which in our case is the edge set of the complete graph $K_n$ on $n$ vertices. A family $\cG$ of winning sets is given in advance, usually containing representatives of a graph theoretic structure, e.g.~all spanning trees or all copies of a fixed graph. Maker wins the game if he occupies all elements of a winning set $G\in \cG$, and Breaker wins otherwise, i.e.~if all the elements of $X$ are claimed and Maker did not fully occupy any $G\in \cG$. 
There is a vast literature on positional games on graphs, we refer the reader to the books \cite{B2008} and \cite{HKSS}. In our setting it is worth mentioning the so-called scoring positional games, where Maker wants to claim as many winning sets as possible, see e.g.~\cite{bdddegop22}.

In Hajnal's triangle game two players also claim unclaimed edges of a complete graph in turns, but this time the graph containing \emph{all} the edges claimed by \emph{both} players should be triangle free. In the original version of the game~\cite{CHR, GP,P}, the player that cannot move, i.e.~who is forced to create a triangle, loses. Afterwards, the extremal version of the game, the \emph{saturation game} was introduced~\cite{BHW, FRS1991, S1992}: the players still have to make sure that no triangles are created in the graph containing all claimed edges, but this time one of the players aims to postpone the game's end as long as possible, i.e.~he tries to maximize \textit{the score of the game}, the total number of edges picked during the game, while the opponent tries to minimize the score. This game has been generalized to arbitrary graphs \cite{CKRW, HKNS, LR, PV}. 

In our Constructor-Blocker games each player builds his own graph, as it is the case in Maker-Breaker games, and in addition Blocker, like Breaker, has no restrictions on his moves. As for Constructor, his graph must remain $F$-free resembling the setting in saturation games.

\medskip 

Note that the following simple observation trivially holds.

\begin{observation}\label{triv1} For two graphs $H,F$ and $n \ge 1$ we have
$$g(n,H,F) \le ex(n,H,F).$$
\end{observation}

In this paper, most of our interest will be directed towards the paths $P_r$ on $r$ vertices and the stars $S_r$ with $r$ leaves. These are  families very much studied in extremal graph theory. The case $H=P_2=S_1=K_2$ corresponds to maximizing/minimizing the number of edges claimed by Constructor. A trivial strategy of Constructor would be to tell Blocker in advance which particular copy of an $F$-free subgraph $G$ of $K_n$ with $ex(n,F)$ edges he will play in, and no matter what Blocker does, Constructor will be able to claim half of the edges of this copy of $G$. Together with Observation \ref{triv1}, we obtain the following.

\begin{observation}\label{triv2} For a graph $F$ and $n \ge 1$ we have
$$\frac{1}{2}ex(n,F)\le g(n,K_2,F) \le ex(n,F).$$
\end{observation}

Our results on the case $F=S_{k+1}$ will include determining $g(n,K_2,S_{k+1})$, and our result will show that its value will be either equal or very close to $ex(n,S_{k+1})$, but we do not think that would be the case for general $F$.

\subsection{Our results}

\paragraph{Star-star games.} First we consider the case $H=S_\ell$ and $F=S_{k+1}$, for some $1\le \ell\le k$. If $\ell\ge 2$, then the score of the game is $\sum_{v\in V(K_n)}\binom{d(v)}{\ell}$, where $d(v)$ is the degree of $v$ in Constructor's graph at the end of the game. If $\ell=1$, then we count the number $\frac{1}{2}\sum_{v\in V(K_n)}d(v)$ of edges.
For comparison, we state the existence of $k$-regular or \textit{almost $k$-regular graphs} (graphs with all but one vertex having degree $k$ and the last one having degree $k-1$) as these graphs contain the maximum number of $S_\ell$'s and thus  determine the value of $ex(n,S_\ell,S_{k+1})$.

\begin{theorem}\label{exstars} For any $k \ge \ell\ge 2$ we have 
    $$ex(n,S_{\ell},S_{k+1})= \begin{cases}
                        \binom{k}{\ell} \cdot n & \textrm{if } nk \textrm{ is even,} \\ 
                        & \\
                        \binom{k}{\ell}(n-1) + \binom{k-1}{\ell} & \textrm{if } nk \textrm{ is odd.}
                         \end{cases} $$ 
If $k\ge \ell=1$, then $ex(n,S_1,S_{k+1})=\lfloor \frac{kn}{2}\rfloor$.
\end{theorem}

Note that with his last claimed edge Blocker can prevent Constructor from building a $k$-regular graph whenever $nk$ is even, so he can achieve that the sum of degrees in Constructor's graph is at most $nk-2$ or $nk-1$, depending on the parity of $nk$. This gives the upper bound for the game score, and we are able to provide a matching lower bound.

\begin{thm}\label{starsvsstars} 



For  
$2\le \ell \le k$  there exists $n_0(k,\ell)$ such that if $n\ge n_0(k,\ell)$ we have 

$$g(n,S_{\ell},S_{k+1})=\begin{cases}
                        \binom{k}{\ell} \cdot (n-2) + 2 \cdot \binom{k-1}{\ell} & \textrm{if } nk \textrm{ is even,} \\ 
                        & \\
                        \binom{k}{\ell}(n-1) + \binom{k-1}{\ell} & \textrm{if } nk \textrm{ is odd.}
                         \end{cases} $$
For $1=\ell\le k$, we have $g(n,S_1,S_{k+1})=\lfloor \frac{nk-1}{2}\rfloor$.

\end{thm}

\paragraph{Tree-star games.} Next we consider the case when $H$ is a tree $T$ and $F=S_{k+1}$. The number of vertices in a graph $T$ is denoted by $|T|$. Clearly, if the maximum degree of $T$ is more than $k$, then $ex(n,T,S_{k+1})=0$. Observe that if the girth of a $k$-regular graph $G$ on $n$ vertices is more than $|T|$ (the radius of $T$ would suffice), then the number of copies of $T$ in $G$ is the maximal that can be achieved by greedily embedding $T$ to $G$. The existence of such graphs is well-known, implying the following extremal result.

\begin{theorem}\label{B} For any $ k \ge 2$ and tree $T$ with maximum degree at most $k$, there exists $n_k(T)$ such that whenever $n\ge n_k(T)$, we have
$$ex(n, T, S_{k+1})=m_{n,k,T},$$ where $m_{n,k,T}$ is the number of $T$'s in an almost $k$-regular graph on $n$ vertices with girth more than $|T|$.
\end{theorem}

We can prove that the score of the game is not far from $ex(n, T, S_{k+1})$.

\begin{thm}\label{pathsvsstars} For any $ k \ge 2$ and tree $T$ with maximum degree at most $k$,  we have
$$g(n,T,S_{k+1})=ex(n,T,S_{k+1}) - O_{k,T}(1).$$ 

\end{thm}

\paragraph{Path-path games.} We consider two games with both $H$ and $F$ being paths.

One can find generalized Tur\'an results about paths in the article of Gy\H{o}ri, Salia, Tompkins and Zamora~\cite{GySTZ19}.

\begin{theorem}[\cite{GySTZ19}, Theorem 14, Remark 3\footnote{In \cite{GySTZ19}, the value of $ex(n,P_4,P_5)$ is incorrectly stated to be $\left \lfloor \frac{n-1}{2} \right \rfloor \left \lceil \frac{n-1}{2} \right \rceil$.}] There exists an integer $n_0$ such that for all $n \ge n_0$ we have
\[
ex(n,P_3,P_4)=\binom{n-1}{2}.
\]

For every integer $n\ge n_0$ we have
\[
ex(n,P_4,P_5)= \left \lfloor \frac{n-2}{2} \right \rfloor \left \lceil \frac{n-2}{2} \right \rceil.
\]
\end{theorem}

Let
\[
B(n):=\binom{\left\lfloor \frac{n-2}{2} \right\rfloor}{2} + \binom{\left\lceil \frac{n-2}{2} \right\rceil}{2}.
\]
In the case $H=P_3, F=P_4$ we can determine the exact game score.

\begin{thm}\label{p3p4}
There exists an integer $n_0$ such that for every $n\geq n_0$, we have
\[
g(n,P_3,P_4)=B(n).
\]
\end{thm}

When $H=P_4$ and $F=P_5$ we provide the following bounds.
\begin{thm}\label{p4p5} We have
\[ 
\frac{8}{49}n^2 - o(n^2) \leq g(n,P_4,P_5)\leq \frac{4}{23} n^2 + o(n^2).
\]
\end{thm}

\paragraph{Triangle-path game.} Finally, we study the first instance of the case $H=K_3$, and $F=P_k$. Observe that all components of a $P_4$-free graph that contain a triangle are triangles, so Blocker can easily prevent Constructor building any triangles, i.e.~$g(n,K_3,P_4)=0.$ Therefore we consider the case $H=K_3$, $F=P_5$. Luo proved the following.

\begin{theorem}[\cite{L}, Corollary 1.7] We have
$$ex(n,K_3,P_5)=n-O(1).$$
\end{theorem} 

The extremal graph consists of pairwise vertex-disjoint $K_4$'s. Our next result shows that Constructor can build one out of the four possible  $K_3$'s in each $K_4$.

\begin{thm}\label{k3p5} We have
$$\displaystyle g(n,K_3,P_5)= \frac{n}{4} - o(n).$$
\end{thm}

The rest of the paper is organized as follows. In Section~\ref{s:psss} we deal with the tree-star games and the star-star games, proving Theorem~\ref{starsvsstars} and Theorem~\ref{pathsvsstars}. Section~\ref{s:pp} is devoted to the path-path games, where we prove Theorem~\ref{p3p4} and Theorem~\ref{p4p5}. Finally, the proof of Theorem~\ref{k3p5} on the triangle-path game is in Section~\ref{s:tp}.

\medskip

\textbf{Notation.}

In a Constructor-Blocker game, given $i \ge 1$, we denote by $G_i[\cC]$ ($G_i[\cB]$) the graph of the first $i$ edges picked by Constructor (Blocker, resp.) and we denote by $G_{\textrm{end}}[\cC]$ the graph formed by the edges picked by Constructor at the end of the game.

Given a graph $G$, we denote by $\overline{d}(G)$  the average degree of that graph.

\section{Tree-star games and star-star games} \label{s:psss}

In this section, we prove Theorem~\ref{starsvsstars} and Theorem~\ref{pathsvsstars}. The two proofs use the same ideas and the calculations are also similar, therefore we introduce a general framework. The number of $S_\ell$'s in a graph $G$ is $\sum_{v\in V(G)}\binom{d_G(v)}{\ell}$ if $\ell\ge 2$, and $\frac{1}{2}\sum_{v\in V(K_n)}d(v)$ if $\ell=1$, so to maximize this quantity, Constructor, who is not allowed to have degree $k+1$ in his graph, should build a graph that is "as close to being $k$-regular as possible". When Constructor wants to maximize the number of copies of $T$ for some tree $T$ on $t$ vertices, then based on Theorem \ref{B}, in addition to building an almost $k$-regular graph, he has to make sure that the girth of his graph is larger than the radius of $T$.

For a non-negative integer $C$, let us introduce an auxiliary game that we refer to as the \textit{$S_{k+1}$-free $C$-bounded symmetric forbidden neighborhoods game} ($S_{k+1}$-free $C$-BSFN game, for short). In this game, most things stay the same, Constructor is still not allowed to create an $S_{k+1}$. But on top of this, in every move he is shown a family $\{F_v:v\in K_n \}$ of forbidden neighborhoods with the property that $u\in F_v$ if and only if $v\in F_u$ and all $F_u$ have size at most $C$, and he is not allowed to pick an edge $uv$ with $u\in F_v$. These $F_v$'s may change from move to move, but their size is never more than $C$. The ordinary $S_{k+1}$-free game is equivalent to the case $C=0$. 

With the help of this game we will later help Constructor to maintain high girth, by forbidding edges joining two vertices at distance at most $t-1$. As all degrees in Constructor's graph are at most $k$ throughout the game, $C=k(k-1)^{t-2}$ will be an adequate choice for that game.

Some of our auxiliary lemmas will be stated in this context as well. A graph $G=(V,E)$ together with sets $F_v\subseteq V$ with $|F_v|\le C$ for all $v\in V$ will be called a \textit{graph with $C$-BSFN.}

Our general result is as follows.

\begin{thm}\label{boxthm} 
For any $k \ge 1$, $C\ge 0$ there exists $n_0=n_0(k,C)$ such that if $n\ge n_0$, then in the $S_{k+1}$-free $C$-BSFN
game on a vertex set of size $n$, Constructor can build a graph of minimum degree at least $k-1$ such that $|\{v:d_{G[\cC]}(v)=k-1\}|\le 2+4C$. Moreover, if $C=0$ and $nk$ is odd, then $|\{v:d_{G[\cC]}(v)=k-1\}|\le 1$. 
\end{thm}

For the proof of Theorem \ref{boxthm}, we will need four technical lemmata. Before introducing these auxiliary statements, let us briefly summarize Constructor's strategy so that the Reader should have a better understanding on why these lemmas will turn out to be helpful. It is easy to come up with a strategy for Constructor to build an almost perfect matching playing on $K_n$ even if Blocker starts the game. Therefore Constructor's strategy in the C-BSFN game will consist of two parts: in the first part, Constructor creates a graph with most degrees being $k-1$ and the rest of them $k$. Then in the second part, Constructor tries to build an almost perfect matching on the vertices of degree $k-1$. Of course, after the first part, players do not play on $K_n$ any more, as some edges are already taken, so Constructor needs to come up with a modified strategy for some simple scenarios. This is done in Lemma \ref{tuchlummu3} and Lemma \ref{tuchlummu4}, when the assumption on the graph of occupied edges is either that at least half of the vertices are isolated or that all vertices are of bounded degree. So when building his graph with minimum degree $k-1$ (during the first part of his strategy), Constructor must make sure that every vertex of degree $k-1$ in his graph, should be adjacent to many free edges. This is done in Lemma \ref{firststep}. Then a technical analysis will be presented to show how Constructor can achieve from this initial graph that the total graph should satisfy the assumptions of Lemma \ref{tuchlummu3} and Lemma \ref{tuchlummu4}.

\begin{lemma}\label{firststep}

For any integers $k\ge 1$, $C\ge 0$ and $\varepsilon>0$ there exists  $n_0(k,C,\varepsilon)$ such that for $n \ge n_0(k,C,\varepsilon)$ Constructor can play on $K_n$ in the $S_{k+1}$-free $C$-BSFN game such that after some round $t \ge 1$ we have the following:
\begin{enumerate}

    \item $d_{G_t[\cC]}(v) \in \{k-1,k\}$ for all $ v \in V(K_n)$,

    \item if $d_{G_t[\cC]}(v)=k-1$, then $d_{G_t[\cC] \cup G_t[\cB]}(v) \le \varepsilon n$,

    \item $|\{ v \in V(K_n): d_{G_t[\cC]}(v)=k \}| \le  \varepsilon n$.
\end{enumerate}

\end{lemma}

\begin{proof}
We say that a vertex $v$ is \textit{dangerous} if $d_{G_i[\cC] \cup G_i[\cB]}(v) \ge \varepsilon n/2$. Constructor's strategy is as follows: if after round $i$ there exists a dangerous vertex $v$ with $d_{G_i[\cC]}(v)<k$, then Constructor considers an arbitrary such $v$ and picks an available edge $uv$ with $u\notin F_v$ such that among such vertices $u$, the Constructor-degree $d_{G_i[\cC]}(u)$ is minimum. If after round $i$ all dangerous vertices $w$ have $d_{G_i[\cC]}(w) = k$, then Constructor considers a vertex $v$ with minimum $d_{G_i[\cC]}(v)$ and picks an available $uv$ with $u\notin F_v$ such that among such vertices $u$, the Constructor-degree $d_{G_i[\cC]}(u)$ is minimum.

We claim that for every $j=0,1,\dots,k-1$, there exists $i_j$ such that 

\medskip 
(1) $d_{G_{i_j}[\cC]}(v) \in \{j,j+1\}$ for all non-dangerous $ v \in V(K_n)$, 

\smallskip 
(2) if $d_{G_{i_j}[\cC]}(v)=j$, then $d_{G_{i_j}[\cC] \cup G_{i_j}[\cB]}(v) \le \varepsilon n$, 

\smallskip 
(3) $|\{ v \in V(K_n): d_{G_{i_j}[\cC]}(v)=j+1 \}| \le  \varepsilon n$, and 

\smallskip 
(4) $d_{G_{i_j}[\cC]}(v)=k$ for all dangerous vertices $v\in V(K_n)$. 

\medskip 
Once this claim is proved, the statement of the lemma is the case $j=k-1$. We will prove the claim by induction on $j$, with the case $j=0$ trivially true at the beginning of the game.

 First observe that, as the number of edges played during the game is at most $kn$, the number of vertices that become dangerous during the game is at most a constant $D=D(k,\varepsilon)$. Next, we claim that after any round $i$ for a vertex $v$ with $d_{G_{i}[\cC]}(v)<k$, the number of vertices $u$ for which either $uv$ is already picked or $u\in F_v$ is at most $3\varepsilon n/4$ using the trivial bound $d_{G_{i}[\cC]\cup G_i[\cB]}(v)+C$ for this quantity. This is certainly true for non-dangerous vertices, as for such $v$ by definition we have $d_{G_{i}[\cC]\cup G_i[\cB]}(v)+C\le \varepsilon n/2+C\le 3\varepsilon n/4$ for $n$ large enough. At the moment a vertex $v$ becomes dangerous, its total degree is $\varepsilon n/2$. By the strategy of Constructor, until $v$ reaches Constructor degree $k$, there can be at most $kD$ turns, and thus $d_{G_{i}[\cC]\cup G_i[\cB]}(v)\le \varepsilon n/2+2kD$ and so $d_{G_{i}[\cC]\cup G_i[\cB]}(v)+C\le \varepsilon n/2+2kD+C\le 3\varepsilon n/4$ for $n$ large enough. 
 
 So assume our claim above holds for $j-1$. Then, because of the previous observation, as long as $|\{v:d_{G_i[\cC]}(v)=j-1\}|>3\varepsilon n/4$, Constructor after choosing $v$ will always pick an edge $uv$ with $d_{G_i[\cC]}(u)=j-1$, and so all non-dangerous vertices will have Constructor degree $j-1$ or $j$. This also implies that when the game reaches $|\{v:d_{G_i[\cC]}(v)=j-1\}|\le 3\varepsilon n/4$, we will have $|\{v:d_{G_i[\cC]}(v)=j\}|>(1-\varepsilon) n$, and thus during this phase the "$u$-vertex" of the edge picked by Constructor will have $d_{G_i[\cC]}(u)=j-1$ or $j$. By our strategy, the $v$-vertex is either dangerous or has $d_{G_i[\cC]}(u)=j-1$. Therefore until all $j-1$ Constructor-degree vertices and all dangerous vertices with Constructor-degree less than $k$ are eliminated, at most $3\varepsilon n/4+kD<\varepsilon n$ vertices of Constructor-degree $j+1$ are created. This finishes the induction step and thus the proof of the lemma.
\end{proof}

\begin{lemma}\label{tuchlummu2}
Let $G$ be a graph with $C$-BSFN and with average degree $d$, maximum degree $\Delta$. Also suppose $|V(G)| \ge 1+2d+2C+ \Delta$. Then there exist two non-adjacent vertices $x$ and $y$ in $G$ such that $d_G(x)+d_G(y) \ge 2d$ and $x\notin F_y$. 
\end{lemma}

\begin{proof}
Suppose for a contradiction that $G$ is a counterexample. Let $Z$ denote the set of vertices with degree at least $d$, i.e.\ $Z:=\{v : d_G(v) \ge d\}$. If we can find two non-adjacent vertices $x,y$ in $Z$ with $x\notin F_y$, then we are done. 

Then we have $\Delta< 2d$ since otherwise an $x$ of degree at least $2d$ and any non-adjacent $y$ with $y\notin F_x$ (and such $y$ exists as $|V(G)| \ge 1+2d+2C + \Delta$) would show that $G$ is not a counterexample. This argument also proves that $Z\subseteq \{v\}\cup N_G(v)\cup F_v$ holds for all $v\in Z$, in particular $|Z| \le 1 + \Delta +C$.

Now consider the auxiliary bipartite graph $B$ with classes $Z$ and $V(G)\setminus Z$ with $uv\in E(B)$ if and only if $uv \notin E(G)$, $u \in Z, v \notin Z$ and $u\notin F_v$. As $Z\subseteq \{v\}\cup F_v\cup N_G(v)$ and the maximum degree of $G$ is at most $|V(G)| - 1 - 2d -2C$, therefore $$d_B(v)\ge
|V(G)|-1-|F_v|-d_G(v)\ge 2d+C\ge 1+\Delta +C\ge |Z|$$ for any $v\in Z$, where for the last two inequalities, we used the inequalities ($\Delta<2d$ and $|Z|\le C+1+\Delta$) proved in the previous paragraph. So by Hall's condition, there exists a matching $M$ that covers $Z$. If there  exists $xy \in M$ with $d_G(x)+d_G(y)\ge 2d$, then $G$ is not a counterexample. Otherwise $$d=\frac{1}{|V(G)|}\sum_{x}d_G(x)=\frac{1}{|V(G)|}\left(\sum_{xy\in M}(d_G(x)+d_G(y))+\sum_{v\notin \bigcup_{e\in M}e}d_G(v) \right)<d.$$ 
This contradiction completes the proof.
\end{proof}

\begin{lemma}\label{tuchlummu3}
Let $n,C$ be integers with $n\ge 4C\ge 0$, $n\ge 3$ and suppose  $G\subseteq K_n$ is a graph with $V(G)=V(K_n)$ such that $D(G):=\{v\in V(G):d_G(v)>0\}$ is of size at most $\frac{n-C}{2}$. Then in a Blocker-start $C$-SBFN game, Constructor can build a matching in $K_n\setminus G$ with at least $\lfloor\frac{n-4C-1}{2}\rfloor$ edges.
\end{lemma} 

\begin{proof}
We proceed by induction on $n$. If $C\ge 1$ and $n=4C$ or $n=4C+1$, then there is nothing to prove. If $C=0$ and $n=3$ or $n=4$, then $G$ contains at most one edge, so Constructor can take an edge after Blocker's first move.

Suppose $n\ge 4C+2$ and the statement is proved for $n-2$ and $G\subseteq K_n$ is as in the statement of the lemma. We consider two cases, depending on whether the edge taken in Blocker's first move is disjoint with $D(G)$ or not.

 If Blocker plays an edge $uv$ with $u\in D(G)$, then Constructor can take any edge $uv'$ with $v\neq v'$, $v'\notin F_u$. Such $v'$ exists as $|V(K_n)\setminus[ D(G)\cup\{v\}]|\ge \lceil \frac{n+C}{2}\rceil-1>C$. Then $G\setminus\{u,v'\}$ and $K_n\setminus\{u,v'\}$ satisfy the induction hypothesis.

Suppose that Blocker plays an edge $uv$ with $u,v\notin D(G)$. If $|D(G)|=\lfloor \frac{n-C}{2}\rfloor$, then $4C\le n$ implies $C<\lfloor\frac{n-C}{2}\rfloor=|D(G)|$, and thus Constructor can spot $v'\in D(G)$ with $v'\notin F_u$. Constructor can play $uv'$, and $G\setminus\{u,v'\}$, $K_n\setminus \{u,v'\}$ and $D(G\setminus \{u,v'\})=D(G)\setminus \{v'\}$ satisfy the induction hypothesis. Finally, if $|D(G)|<\lfloor \frac{n-C}{2}\rfloor$, then $|D(G)|\le \frac{n-2-C}{2}$, and Constructor can play any $uv'$ with $v'\neq v$, $v'\notin F_u$. Then $D(G\setminus \{u,v'\})=D(G)$ is still small enough to satisfy the induction hypothesis for $G\setminus \{u,v'\}, K_n\setminus\{u,v'\}$.
\end{proof}

\begin{lemma}\label{tuchlummu4}
Let $G \subseteq K_n$ with $V(G)=V(K_n)$ be a graph with maximum degree at most $\Delta$. Then in a Blocker-start $C$-BSFN game, Constructor can build a matching $M$ in $K_n\setminus G$ covering at least $n-\Delta-2-C$ vertices. Moreover, if $\Delta\le 1, C=0$ and $n$ is odd, then $M$ covers $n-1$ vertices, while if $\Delta\le 1, C=0$ and $n$ is even, then $M$ covers $n-2$ vertices.
\end{lemma}

\begin{proof}
Suppose first $C>0$ or $\Delta \ge 2$. If $n< \Delta+C+3$, there is nothing to prove. 
Then we proceed by induction on $n$. If $G \subseteq K_n$ has maximum degree at most $\Delta$ and Blocker's first edge is $xy$, then, as $n \ge \Delta + 3 + C$, there exists a $z\in V(K_n) \setminus \{y\}$ with $xz\notin E(G)$, $z\notin F_x$. Constructor can play the edge $xz$, and then $G[V(G)\setminus \{x,z\}]\subseteq K_{n-2}$ has maximum degree at most $\Delta$ and the statement follows.

The case $\Delta=C=0$ is covered by Lemma \ref{tuchlummu3}. If $C=0, \Delta=1$, then $G$ is a partial matching. If $n=2$, then there is nothing to prove. If $n=3$, then $G$ consists of a single edge, so Constructor can claim the third edge in $K_3$ even after Blocker's starting move. For $n\ge 4$, we proceed by induction on $n$: if Blocker picks $xy$, then, as $n\ge 4$ and $\Delta= 1$, there exists $z$ with $xz$ unclaimed, so Constructor can pick $xz$. Then $G[V(G)\setminus \{x,z\}]$ satisfies the inductive hypothesis.
\end{proof}

\begin{proof}[Proof of Theorem \ref{boxthm}] \ 

\textbf{Phase 1.} First by choosing $\varepsilon = \frac{1}{100k^2}$ in Lemma \ref{firststep}, Constructor can build a graph such that after round $t_1$ we have:
\begin{enumerate}

    \item $d_{G_{t_1}[\cC]}(v) \in \{k-1,k\}$ for all $ v \in V(K_n)$,

    \item if $d_{G_{t_1}[\cC]}(v)=k-1$, then $d_{G_{t_1}[\cC] \cup G_{t_1}[\cB]}(v) \le \frac{n}{100k^2}$, and

    \item $|\{ v \in V(K_n): d_{G_{t_1}[\cC]}(v)=k \}| \le \frac{n}{100k^2}$.


\end{enumerate}

\medskip 

\textbf{Phase 2.} The goal of Constructor during Phase 2 (that starts with round $t_1+1$) is to decrease below 2 the average total degree in the graph that is induced by the vertices with Constructor degree $k-1$. To do so Constructor will consider the graph that contains vertices of Constructor degree $k-1$, and pick the vertices $x$ and $y$ to connect relying on Lemma \ref{tuchlummu2}. 

\medskip 

Let us denote by $X_t^{k-1}(\cC)$ 
the set of those vertices whose Constructor degree is $k-1$  after round $t$, i.e.\ $X_t^{k-1}(\cC):=\{v : d_{G_t[\cC]}(v)=k-1\}$. For $x \in X_t^{k-1}(\cC)$ let us denote by $d_t(x)$ the total degree of $x$ in the graph induced by the vertices of $X_t^{k-1}(\cC)$, so in the graph $G_t[X_t^{k-1}(\cC)]$. Let us denote by $\overline{d}_t$ the average total degree and by $\Delta_t$ the maximum total degree of the graph $G_t[X_t^{k-1}(\cC)]$.

\medskip 

Using these notations, after round $t_1$ we have: 

\medskip 

$\bullet$ $|X_{t_1}^{k-1}(\cC)| \ge \left(1-\frac{1}{100k^2}\right)n$, 

\smallskip 

$\bullet$ $\Delta_{t_1} \le \frac{n}{100k^2}$, and 

\smallskip 

$\bullet$ $\overline{d}_{t_1}\le \frac{2\cdot \# \{\text{total number of edges played}\}}{|X^{k-1}_t(\cC)|} \le \frac{2\sum_vd_{G_t[\cC]}(v)}{|X^{k-1}_t(\cC)|}\le \frac{2(k-1+\frac{1}{100k^2})n}{(1-\frac{1}{100k^2})n}\le 2k$. 

\medskip 

In each round $t$ with $t_1 \le t \le t_2$, for some later defined $t_2$ that we will need to be at most $t_1 + \frac{(1-\varepsilon)n(k-1)}{2k}$,  Constructor applies Lemma~\ref{tuchlummu2} with $d=\overline{d}_t$, $\Delta=\Delta_t$, 
finds vertices $x$ and $y$ and connects them with an edge. If there are more possibilities (of pairs of vertices) to pick, then he prioritizes an edge that contains a vertex with maximum total degree among vertices that are adjacent to edges $xy$ satisfying the statement of Lemma~\ref{tuchlummu2}.


Let us suppose that Constructor connects the picked vertices (and so these vertices will have Constructor degree $k$). The analysis how the average total degree changes is 
the next claim. Consider now the maximum degree. Blocker can increase the degree of one or two vertices $x$ and $y$ of maximum total degree. If this increased maximum degree is still smaller than $2\overline{d}_t$, then no matter what Constructor plays, we will have $\Delta_{t+1}\le 2\overline{d}_t$. If this increased maximum degree is larger than $2\overline{d}_t$, then by the prioritization rule of Constructor's strategy, Constructor will pick an edge adjacent to $x$. Then $x$ and the edge $xy$ will be eliminated, and the degree of $y$ will be again at most $\Delta_t$. We obtained the following claim.

\begin{claim}\label{Deltad} For all $t$ with $t_1 \le t \le t_2$ (where  $t_2 \le t_1 + \frac{nk}{2(k+1)})$, we have

            \begin{equation}\label{eq_Dt} 
                \Delta_{t+1} \le \max\{\Delta_{t}, 2\overline{d}_t\}\le \frac{n}{100k^2}, \textrm{ and}
            \end{equation}  
            \begin{equation}\label{eq_dt}
                \overline{d}_{t+1} \le \frac{\overline{d}_t|X_t^{k-1}(\cC)|-2(d_t(x)+d_t(y)) + 2}{|X_t^{k-1}(\cC)|-2} \le \frac{\overline{d}_t \cdot (|X_t^{k-1}(\cC)|-4) +2}{|X_t^{k-1}(\cC)|-2}.
            \end{equation}

\end{claim}

For all $t$ with $t_1 \le t \le t_1 + \frac{nk}{2(k+1)}$, the number of vertices in $X_t^{k-1}(\cC)$ is at least $n(1-\frac{1}{100k^2}-\frac{k}{k+1})\ge \frac{n}{2(k+1)}\ge 1+2\overline{d}_t+2c+2\Delta_t$, so the assumptions of Lemma \ref{tuchlummu2} stay valid for these values of $t$.
Therefore as long as we choose $t_2$ with $t_2 \le t_1 + \frac{nk}{2(k+1)}$, Constructor is able to follow the given strategy.

\medskip 

Note that if $\overline{d}_t \ge 2$, then by (\ref{eq_dt}) we have $$\overline{d}_{t+1} \le \overline{d}_t \cdot \frac{|X_t^{k-1}(\cC)|-3} {|X_t^{k-1}(\cC)|-2},$$
so if Constructor follows the given strategy in Phase 2 during $s$ rounds after $t_1$, then either at some round we have $\overline{d}_t<2$ or we have

\begin{equation}\label{eq_telescope}
\overline{d}_{t_1+s} \le \overline{d}_{t_1} \cdot  \frac{|X_{t_1}^{k-1}(\cC)|-3} {|X_{t_1}^{k-1}(\cC)|-2}\cdot\frac{|X_{t_1}^{k-1}(\cC)|-5} {|X_{t_1}^{k-1}(\cC)|-4} \cdots \frac{|X_{t_1}^{k-1}(\cC)|-(2s+1)} {|X_{t_1}^{k-1}(\cC)|-2s} \le \overline{d}_{t_1} \cdot \frac{|X_{t_1}^{k-1}(\cC)|-(2s+1)}{|X_{t_1}^{k-1}(\cC)|-2}.
\end{equation}

As $\overline{d}_{t_1}\le 2k$, and the fraction in the right hand side of (\ref{eq_telescope}) is at most $\frac{n-(2s+1)}{(1-\frac{1}{100k^2})n-2}$, by putting $s_1= \frac{nk}{2(k+1)}$, we have $\overline{d}_{t_1+s_1}<2$, which means that we will have a first round $t_1+s \le t_1 + \frac{nk}{2(k+1)}$ for which $\overline{d}_{t_1+s} < 2$. Let us define $t_2$ as $t_1+s$. Note that at the end of Phase~2 we have $|X^{k-1}_{t_1+s}| \ge (\frac{1}{k+1}-\frac{1}{100k^2})n$. 

\bigskip 

\textbf{Phase 3.} In this phase, Constructor's goal is to achieve an even sparser induced subgraph, obtaining either $\overline{d}_t < \frac{1}{2}$ or the graph contains just a matching. Then he will be able to finish the game with either Lemma \ref{tuchlummu3} or Lemma \ref{tuchlummu4}.

\medskip 

More precisely, in Phase 3 Constructor first does a similar thing as in Phase 2. He tries to pick 2 vertices $x,y \in X_t^{k-1}(\cC)$ with $d_t(x) + d_t(y) \ge 4$ and $y \not \in F_x$ and connect them.  

$\bullet$ If he can find such pairs till the average total degree in $G_t[X_t^{k-1}(\cC)]$ will be less than $\frac{1}{2}-\delta$ with some $\delta$, then he continues as in Case 1 below, or 

$\bullet$ if at some round he can not find 2 vertices $x,y \in X_t^{k-1}(\cC)$ with $d_t(x) + d_t(y) \ge 4$ and $y \not \in F_x$, then he continues either as in Case 2, Case 3 or Case 4 below. 

\smallskip 

Note that if Constructor can pick 2 vertices $x,y \in X_t^{k-1}(\cC)$ with $d_t(x) + d_t(y) \ge 4$ and $y \not \in F_x$ then in (\ref{eq_dt}) we still have $$\overline{d}_{t+1} \le \overline{d}_t \cdot \frac{|X_t^{k-1}(\cC)|-3} {|X_t^{k-1}(\cC)|-2}.$$

\medskip 

\textbf{Case 1.} For all $t$ with $t_2 \le t \le t_3 \le  t_2 + \frac{n}{2(k+2)}$ we have 2 vertices $x,y \in X_t^{k-1}(\cC)$ with $d_t(x) + d_t(y) \ge 4$ and $y \not \in F_x$. 

Then using (\ref{eq_telescope}) we have $\overline{d}_{t_2 + s} < \frac{1}{2}- \delta$ with some $s \le \frac{n}{2(k+2)}$ and $\delta > 0$. It means that the number of vertices that are touched by some edge is less than $|X^{k-1}_{t_2+s}(\cC)|(\frac{1}{2}-\delta)$. As $|X^{k-1}_{t_2+s}(\cC)|(\frac{1}{2}-\delta) < \frac{|X^{k-1}_{t_2+s}(\cC)|-C}{2},$ and $|X^{k-1}_{t_2+s}(\cC)|$ is linear in $n$, we can apply  Lemma \ref{tuchlummu3} after round $t_2 + s$ and we are done.

\bigskip

Note that if we are not in Case 1, for the graph $G_t[X_t^{k-1}(\cC)]$, we have $\Delta_t \le 3$. Indeed,  we know that $\Delta_t \le \varepsilon n$ by (\ref{eq_Dt}), and even if we consider $t=t_2+\frac{n}{2(k+2)}$, the size of the set of vertices of Constructor degree $k-1$ is at least $(\frac{1}{(k+1)(k+2)}-\frac{1}{100k^2})n$, and so as $\Delta_t + 1 + C < |X^{k-1}_{t}(\cC)|$, that means there is a vertex $y \in X_t^{k-1}(\cC)$ such that $y \not \in F_x$, which means we are in Case 1, a contradiction. Therefore Lemma \ref{tuchlummu4} finishes the proof unless $C=0$. So from now on we will assume $C=0$ and we need to deal with the cases $\Delta_t\in\{0,1,2,3\}$.

\bigskip

\textbf{Case 2.} $\Delta_t\le 1$

The moreover part of Lemma \ref{tuchlummu4} finishes the proof.

\bigskip 

\textbf{Case 3.} $\Delta_t=2$

\medskip 


Note that as we are not in Case 1, we cannot have two non-adjacent vertices of degree 2, so we have 1, 2 or 3 vertices of degree two and they form a clique. If $x$ is one such vertex, then Constructor connects $x$ to any available vertex. As degree 2 vertices formed a clique of size at most 3, $G_{t+1}[X^{k-1}_{t+1}(\cC)]$ is a matching and Constructor can apply the moreover part of Lemma \ref{tuchlummu4} and we are done.








%
%

\medskip 

\textbf{Case 4.} $\Delta_t=3$

Let $x$ be a vertex in $X_t^{k-1}(\cC)$ with $d_t(x)=3$. Observe that the vertices in $X_t^{k-1}(\cC)$ that are not adjacent to $x$ form an independent set in $G_t[X_t^{k-1}(\cC)]$ (otherwise we would be in Case 1), so we can apply Lemma \ref{tuchlummu3}, as $4 < \frac{|X_t^{k-1}(\cC)|}{2}$.
\end{proof}

\begin{proof}[Proof of Theorem \ref{starsvsstars}] Constructor follows his strategy in a $S_{k+1}$-free $0$-BSFN
game provided by Theorem \ref{boxthm}. So he can build a graph $G_{\textrm{end}}[\cC]$ of minimum degree $k-1$ with $|\{v:d_{G_{\textrm{end}}[\cC]}(v)=k-1\}|\le 2$. Counting the $S_{\ell}$'s in $G_{\textrm{end}}[\cC]$ confirms the statement of  Theorem \ref{starsvsstars}.
\end{proof}

\begin{proof}[Proof of Theorem \ref{pathsvsstars}]  Constructor follows his strategy in a $S_{k+1}$-free $C$-BSFN
game with $C=k(k-1)^{|T|-2}$ provided by Theorem \ref{boxthm}. At any round, and for any vertex $v$ the forbidden neighborhood $F_v$ consists of all vertices at distance at most $|T|-1$ from $v$. So Constructor can build a graph $G_{\textrm{end}}[\cC]$ of minimum degree $k-1$ with $|\{v:d_{G_{\textrm{end}}[\cC]}(v)=k-1\}|\le 2+4C$ and if $C=0$ and $nk$ is odd, then $|\{v:d_{G_{\textrm{end}}[\cC]}(v)=k-1\}|=1$. Counting the $T$'s in $G_{\textrm{end}}[\cC]$ confirms the statement of Theorem \ref{pathsvsstars}.
\end{proof}


\section{Path-path games} \label{s:pp}

In this section, we prove Theorem~\ref{p3p4} and Theorem~\ref{p4p5}.

\medskip



\begin{proof}[Proof of Theorem~\ref{p3p4}]
Let us observe that any $P_4$-free graph must be a disjoint union of vertices, stars and triangles.

We start by exhibiting a strategy for Constructor that will ensure that at the end of the game his graph consists of two disjoint stars spanning all $n$ vertices. To do that, before the game starts he fixes two vertices, $v_1$ and $v_2$, and then he follows a simple pairing strategy -- whenever Blocker claims an edge $xv_i$, for some vertex $x$ and $i\in\{1,2\}$, Constructor responds by claiming $xv_{3-i}$. If Blocker claims an edge disjoint from $\{v_1, v_2\}$ or an edge $xv_2$ such that $xv_1$ is already claimed by Constructor, and also in the very first move of the game, Constructor picks a vertex $y$ that is isolated in his graph and claims the edge $yv_1$.

Following this strategy Constructor clearly ends up with two disjoint stars on $n$ vertices, $S_{a}$ and $S_{n-a-2}$, for some integer $a$. The number of $P_3$'s in his graph is $\binom{a}{2}+\binom{n-a-2}{2}$, which is minimized when $a$ and $n-a-2$ differ by at most 1, implying $g(n,P_3,P_4)\geq B(n)$.

\medskip

Next, we analyse Blocker's prospects in this game. The following strategy will be referred to as the \emph{Basic Strategy} of Blocker. After each move of Constructor, we locate the connected component $C$ in Constructor's graph containing that move. If $C$ is an isolated edge or a triangle, Blocker responds by claiming an arbitrary edge. If $C$ is a star centered at $v$ with at least two leaves, Blocker responds by claiming an edge incident to $v$, 
if such an edge is available; otherwise, he claims an arbitrary edge.

From the beginning of the game, Blocker will follow the Basic Strategy. He suspends it at most once, when Constructor's graph consists of nontrivial components, $S_1$ and $S_2$, each of which is a star with at least three leaves. For $S_1$ and $S_2$, we denote the center vertices by $v_1$ and $v_2$, respectively, and the number of leaves by $a_1$ and $a_2$, respectively. W.l.o.g.~let us assume that $a_1\geq a_2$. We distinguish two cases.
\begin{itemize}
    \item[(i)] $a_1-a_2 > 3$
    
    Blocker locates a vertex $x$ such that the edge $xv_1$ is already claimed by him and the edge $xv_2$ is unclaimed, and claims $xv_2$. After that, he gets back to following the Basic Strategy to the end of the game.
    
    \item[(ii)] $a_1-a_2 \leq 3$
    
    Blocker repeatedly locates the larger (at that point) of Constructor's two stars, breaking ties arbitrarily, with its center at $v_\ell$, with $\ell\in\{1,2\}$, and claims the edge between $v_{\ell}$ and a vertex $y$ that is isolated in Constructor's graph. If Constructor does not claim the edge $yv_{3-\ell}$ in his following move, Blocker claims that edge and gets back to following the Basic Strategy to the end of the game. Otherwise, he keeps playing according to the strategy described in case (ii).
\end{itemize}
If at any point in the game the strategy calls for Blocker to claim an edge that he claimed earlier, when he was to ``claim an arbitrary edge'', he just claims a new arbitrary edge and continues.

Let us show that Blocker can always follow this strategy. Once the strategy in case (i) is activated, Blocker surely claimed at least $a_1-2$ edges incident to $v_1$ since up to that point he followed the Basic Strategy. As the number of vertices of $S_2$ is $a_2 + 1 < a_1-2$, there must exist a suitable vertex $x$ for Blocker's move in (i). In case (ii), Blocker can clearly play as long as there are isolated vertices in Constructor's graph.

We move on to analysing the Constructor's graph at the end of the game. 

\begin{claim} \label{claim:star}
At the end of the game, Constructor will not have a star with more than $n/2 + 2$ leaves. \end{claim}

\begin{proof}[Proof of Claim]
While Blocker follows the Basic Strategy, at any round Constructor adds an edge to a star Blocker adds his edge to the same center vertex. In every star Constructor has one edge head start when he claims the very first edge, as Blocker does not respond in the same way to isolated edges. Furthermore, due to parity, when Constructor claims the last edge of a star Blocker may not have an available unclaimed edge to respond.

If case (i) is activated, then Constructor plays exactly one move to which Blocker does not respond by the Basic Strategy. Altogether, the number of edges of Constructor incident to a vertex cannot be larger than the number of Blocker's edges incident to that same vertex by more than 3.

If case (ii) is activated, both players claim edges incident to either $v_1$ or $v_2$. As initially $a_1-a_2 \leq 3$, we have that, for both $v_1$ and $v_2$, the number of Constructor's edges and the number of Blocker's edges played incident to that vertex within case (ii) may differ by at most 3. Altogether, the number of edges of Constructor incident to a vertex cannot be larger than the number of Blocker's edges incident to that same vertex by more than 6.

Having in mind that the total number of edges incident to a vertex is $n-1$, the assertion of the claim readily follows.
\end{proof}

Suppose that at the end of the game Constructor has $k$ triangle components, and $t$ star components with, respectively, $s_1, s_2, \dots, s_t$ leaves. We have $3k + \sum_{i=1}^t (s_i+1) \leq n$, and Claim~\ref{claim:star} implies $s_i \leq n/2 +2$, for all $i$. The number of $P_3$'s Constructor created is $N := 3k + \sum_{i=1}^t \binom{s_i}{2}$. 

Let us first analyze the case $k=0$, $t=2$ and $s_1+s_2=n-2$. As the game ends with two stars of Constructor, one of the cases (i), (ii) was activated during the game. Activating case (i) would result in a vertex that is in neither of the two stars, which is in contradiction with $s_1+s_2=n-2$. Similarly, activating case (ii) and afterwards getting back to the Basic Strategy results in a vertex that is in neither of the two stars, which is again not possible. Hence, the only remaining possibility is that Blocker activated case (ii) and stayed in it to the end of the game. But that results in stars which are balanced, i.e.~$|s_1-s_2|\leq 1$. Note that in this case we have $N=B(n)$. Observe that we use the condition $n\ge n_0$ so that Blocker has space to force Constructor to balance his stars with his first two moves played in case (ii).

Our aim is to show that all other possibilities for Constructor's graph at the end of the game, as long as it satisfies the restrictions listed above, result in smaller value of the function $N$. 

If $k\geq 2$, we can replace all triangles with one new star thus increasing $N$. If $k=1$, then removing the triangle and adding its three vertices as new leaves of a star increases $N$. Hence, to maximize $N$ we must have $k=0$.

If $k=0$, $t=2$ and $s_1+s_2<n-2$, Claim~\ref{claim:star} upper bounds both $s_1$ and $s_2$ and in this case the maximal value of $N$ is obtained for $s_1 = \left\lceil n/2+2\right\rceil$ and $s_2=n-3 -s_1$. It is easy to check that such $N$ is less than $B(n)$, for $n$ large enough.

Further on, if $t>2$, then by removing the smallest of the stars and repeatedly attaching its vertices, one by one, to the smallest (at that point) of the remaining stars we increase $N$. Therefore, $t=2$ is required to maximize $N$, as Claim~\ref{claim:star} eliminates the option $t=1$.

Hence, if Blocker follows the described strategy Constructor cannot claim more than $B(n)$ many $P_3$'s, implying $g(n,P_3,P_4)\leq B(n)$.
\end{proof}

\medskip

\begin{proof}[Proof of Theorem~\ref{p4p5}]
Throughout the game, let $\cC$ denote Constructor's graph. We first study Constructor's prospects, exhibiting a strategy enabling him to claim a large double star, and analyzing that strategy to obtain the lower bound. As we are aiming for an asymptotic result, we will omit floors and ceilings for ease of presentation.

In his first move, Constructor claims an edge $v_1v_2$. Then, he proceeds to Stage~1, following a pairing strategy: if Blocker claims the edge $xv_i$, for $i\in\{1,2\}$ and some vertex $x$ isolated in $\cC$, Constructor responds by claiming $xv_{3-i}$, and if Blocker claims an edge disjoint with $\{v_1,v_2\}$, Constructor claims $xv_1$ for some vertex $x$ isolated in $\cC$. Stage~1 ends when one of the vertices $v_1, v_2$ has $3(n-2)/7$ hanging edges in $\cC$, assume w.l.o.g.~that that vertex is $v_1$.

In Stage~2, Constructor repeatedly claims an edge between $v_2$ and a vertex isolated in $\cC$, for as long as such edges are available, and then proceeds to Stage~3.

In Stage~3, Constructor repeatedly claims an edge between $v_1$ and a vertex isolated in $\cC$, for as long as such edges are available. Once this stage ends Constructor is done building his double star and for the remainder of the game he claims arbitrary free edges disjoint from the double star.

Constructor can clearly follow the described strategy, and at the end of the game, one of the components in $\cC$ will be a double star centered at $v_1$ and $v_2$. It remains to estimate the number of $P_4$'s guaranteed to be found in that double star.

After Stage~1 the number of leaves adjacent to $v_1$ in $\cC$ is $3(n-2)/7$, the number of leaves adjacent to $v_2$ in $\cC$ is between $0$ and $3(n-2)/7 - 1$, and at that point Blocker does not have a claimed edge between any vertex isolated in $\cC$ and either $v_1$ or $v_2$.

When Stage~2 finishes, the vertices still isolated in $\cC$ can be classified into two sets: $U_1$, containing those vertices $x$ for which Blocker claimed $xv_2$ but $xv_1$ is unclaimed, and $U_2$, containing those vertices $y$ for which Blocker claimed both $yv_1$ and $yv_2$. Let us denote by $U_0$ the leaves of Constructor's star centered at $v_2$ at the end of Stage~2, and let $u_i:=|U_i|$, $i\in\{0,1,2\}$. Note that 
\begin{equation} \label{e:1}
u_0+u_1+u_2 = \frac{4(n-2)}{7}.
\end{equation}
Furthermore, as Blocker started claiming edges between $U_1 \cup U_2$ and $\{v_1, v_2 \}$ only in Stage~2, and he claimed one for every vertex in $U_1$ and two for every vertex in $U_2$, we have 
\begin{equation} \label{e:2}
u_0\geq u_1+ 2u_2.
\end{equation}

In Stage~3, Constructor will claim at least half of the edges between $v_1$ and $U_1$, and after Stage~3 the number of $P_4$'s in his double star will be at least $\left( 3(n-2)/7 + u_1/2 \right)\cdot u_0$. Expressing $u_2$ from~(\ref{e:1}) and plugging it into~(\ref{e:2}) gives a lower bound on $u_0$, implying that the number of $P_4$'s that Constructor created is at least
\begin{equation*}
    \left( \frac{3(n-2)}{7} + \frac{u_1}{2} \right)\cdot \frac13 \left( \frac{8(n-2)}{7} - u_1 \right) = \frac16 \left( \frac{6(n-2)}{7} + u_1\right) \cdot \left( \frac{8(n-2)}{7} - u_1 \right)\geq \frac{8(n-2)^2}{49}.
\end{equation*}
To get the last inequality we observe the expression preceding it as a quadratic function of $u_1$. This function has its maximum in $u_1=(n-2)/7$, and we want to minimize it.
Knowing that $0\leq u_1\leq 2(n-2)/7$, which is a direct consequence of~(\ref{e:1}) and~(\ref{e:2}), the minimum of the above quadratic function when restricted to that interval is achieved when $u_1=0$.

\medskip

Next, we describe and analyze a strategy for Blocker. Throughout the game we keep track of all nontrivial connected components in $\cC$, and after every move of Constructor, Blocker responds by claiming an edge incident to the component $X$ containing that Constructor's move. Depending on what $X$ is, we distinguish several cases.
\begin{itemize}
\item[(i)] If $X$ is a star centered at $a$, Blocker claims an edge between $a$ and $b$ such that $ab$ is unclaimed, prioritizing vertices $b$ that are centers of other Constructor's star components, sorted by size in decreasing order; if there is no available unclaimed edge incident to $a$, Blocker plays arbitrarily.

\item[(ii)] If $X$ is a double star centered at $u$ and $v$, where the star at $v$ is not larger than the star at $u$, Blocker claims an unclaimed edge between $v$ and a vertex $w$ isolated in $\cC$, prioritizing vertices $w$ for which $uw$ is already claimed by him. If $uw$ is unclaimed and Constructor does not claim it in his following move \emph{in this component,} Blocker responds by claiming $uw$, thus isolating $w$ from that double star. 

If there is no available unclaimed edge at $v$, Blocker plays arbitrarily.

\item[(iii)] If $X$ is neither a star nor a double star, Blocker plays arbitrarily.
\end{itemize}

We claim that whenever Constructor connects stars $S_1$ and $S_2$ in $\cC$ into a double star by claiming the edge between their centers, at least one of them has at most $2\sqrt{n}$ leaves. Suppose for a contradiction that both stars have at least $2\sqrt{n}$ leaves. Looking back earlier in the game when the first of the two stars, w.l.o.g.~let that be $S_1$, got $\sqrt{n}$ leaves, $S_2$ was still to receive more than $\sqrt{n}$ of Constructor's edges. Each of those edges would be responded by a Blocker's edge at the center of $S_2$, as advised by case~(i) of the strategy. But all that time $S_1$ is among the $\sqrt{n}$ largest stars of Constructor, so Blocker will surely claim the edge connecting the centers of $S_1$ and $S_2$ before Constructor does, a contradiction.

Hence, the first time a double star in $\cC$ is created, the smaller of its stars has at most $2\sqrt{n}$ leaves. 
Let us observe the double star that is the largest double star (by the number of vertices) in Constructor's graph at the end of the game. We denote the number of leaves of its smaller star by $s_1$ and the number of leaves of its larger star by $s_2$, where $s_1\leq s_2$.

Throughout the game, case~(ii) of the strategy implies that after the stars are joined into the double star, for every two new leaves of the smaller Constructor's star, Blocker isolates one vertex from that double star. Indeed, looking at the sequence of moves on that component, whenever Constructor increases the smaller star Blocker responds at that same star, following case~(ii). After that point, Constructor may repeatedly increase the other star a number of times, but the next time he goes back to the first star Blocker will isolate a vertex from that Constructor's component. 

We now count the number of times Constructor added a new leaf to the currently smaller star (i.e.~the star that was smaller of the two at the time the leaf was added), from the moment the double star is first created to the end of the game. Observe that the size of the currently smaller star increases only when Constructor adds a leaf to the currently smaller star, and it increases by exactly one. Therefore, this must have happened at least $s_1 - 2\sqrt{n}$ times.

Hence, every increase of the size of the smaller star by two results in (at least) one isolated vertex. If the number of isolated vertices is denoted by $\ell$, we have $\ell \geq (s_1 - 2\sqrt{n})/2$. As $\ell + s_1 + s_2 \leq n$, we have 
\begin{equation} \label{e:a1}
\frac{3s_1}{2} + s_2 \leq n + \sqrt{n}.    
\end{equation}


Let us first assume that $s_2< 2n/11$. In that case, for any Constructor's double star that has $a$ leaves in one star and $b$ leaves in the other, $a+b<4n/11$ holds. Therefore, the number of $P_4$'s in that double star is
$ ab \leq \left(\frac{a+b}{2}\right)^2 < \frac{n}{11}(a+b)$, so the total number of $P_4$'s in the Constructor's graph is upper bounded by $(n/11)\cdot n < 4n^2/23$ and we are done. Hence, from now on we can assume that $s_2\geq 2n/11$.

In case $s_1< 2n/11$, the number of $P_4$'s in the largest Constructor's double star at the end of the game is $s_1s_2 \leq \frac{2n}{11}\left( s_1 + s_2 -\frac{2n}{11}\right)$. For any other Constructor's double star that has $a$ leaves in one star and $b$ leaves in the other, the number of $P_4$'s is $ab \leq \left(\frac{a+b}{2}\right)^2 < \frac{2n}{11}(a+b)$,
having in mind that for all Constructor's double stars that are not the largest we have $a+b< n/2<8n/11$. Hence, the total number of $P_4$'s in Constructor's graph is upper bounded by $(2n/11)(9n/11) < 4n^2/23$ and again we are done, so from now on we can assume that $s_1\geq 2n/11$.

If we denote the number of vertices that are not in the largest Constructor's component by $s_3$, we have $s_3 \leq n - s_1 -s_2$, and the number of $P_4$'s on $s_3$ vertices cannot be larger than $(s_3/2)^2$.  If by $S$ we denote the total number of Constructor's $P_4$'s at the end of the game, we have
\begin{equation*}
S\leq s_1s_2  + \left(\frac{s_3}{2}\right)^2 \leq s_1s_2  + \frac14\left(n - s_1 -s_2\right)^2.
\end{equation*}
From~(\ref{e:a1}) we get $s_2 \leq n - 3s_1/2 +o(n)$, and $s_1\geq 2n/11$ guarantees $s_2 \geq n - 9s_1/2$. Within that interval the previous upper bound on $S$ is maximized for $s_2 = n - 3s_1/2+o(n)$, and we get 
\begin{equation*}
S \leq s_1\cdot\left(n - \frac{3s_1}{2} \right) + \frac14\left(n - s_1 -\left(n - \frac{3s_1}{2} \right)\right)^2 +o(n^2) = s_1n - \frac{23s_1^2}{16} +o(n^2).
\end{equation*}
The last expression is maximized for $s_1 = 8n/23$, giving $S\leq 4n^2/23 +o(n^2)$.
\end{proof}

\section{Triangle-path game} \label{s:tp}

\begin{proof}[Proof of Theorem~\ref{k3p5}]
Note that in a $P_5$-free graph, the only way to have more than one triangle in a connected component is to have a component that is a subgraph of $K_4$.

We start by exhibiting a strategy for Blocker:
\begin{itemize}
    \item[(a)] Whenever Constructor claims an edge that is (after he claimed it) in a Constructor's connected component on three vertices, if there is a free edge in that component, Blocker claims it. 
    \item[(b)] Whenever Constructor claims an edge that is (after he claimed it) in a Constructor's connected component on four vertices, if there is a free edge in that component, Blocker claims it.
    \item[(c)] Otherwise, Blocker plays arbitrarily.
\end{itemize}

If Blocker follows this strategy, because of part (a), Constructor will be unable to ever create a triangle in a component on three vertices. Therefore, when he first creates a component on four vertices it will have three of his edges. The part (b) ensures that after that, as long as this component spans four vertices, Constructor can claim at most one more edge in it. This is not enough to ever create two triangles in a component.

To sum up, Constructor will not create an isolated triangle, and in each of his connected components he will have at most one triangle, implying $g(n,K_3,P_5) \leq n/4$.

Next, we describe and analyze a strategy for Constructor. As long as the number of vertices isolated in his graph is not less than $5\sqrt{n}$, Constructor repeatedly goes through the following sequence of steps, building his connected components one by one.
\begin{enumerate}
    \item Among the isolated vertices in his graph, he spots a set of five vertices $X=\{v_1, v_2,\dots,v_5\}$ such that all edges induced on $X$ are unclaimed, and claims $v_1v_2$.
    If Blocker responds by claiming an edge in $X$, w.l.o.g.~let us assume that this edge is incident with $v_5$.
    \item Next, Constructor claims $v_1v_3$. If Blocker does not respond by claiming $v_2v_3$, Constructor claims it in his following move completing a component on three vertices with a triangle.
    \item Otherwise, Constructor claims $v_1v_4$, and in his following move he claims either $v_2v_4$ or $v_3v_4$, completing a component on four vertices with a triangle.
\end{enumerate}

Clearly, as long as Constructor can follow this strategy his graph will consist of connected components on either three or four vertices, each containing a triangle, and the number of the edges he claimed will be less than $n$. It remains to verify that each time when he is to play step 1, he can spot a set $X$ that satisfies the conditions. 

Denote by $L$ the set of vertices isolated in Constructor's graph, and let $\ell:=|L|$. The total number of five vertex subsets from $L$ is $\binom{\ell}{5}$, and each edge Blocker claimed can be contained in at most $\binom{\ell}{3}$ such sets. Therefore, knowing that $\ell > 5 \sqrt{n}$ and that Blocker claimed less than $n$ edges, there must always be a suitable $X$ completely free of Blocker's edges.

Hence, following the above strategy Constructor will create connected components on either three or four vertices, each containing a triangle, until he spans all but at most $5\sqrt{n}$ vertices. This implies $g(n,K_3,P_5) \geq n/4 - 5\sqrt{n}/4 =n/4 - o(n)$.
\end{proof}

\section{Conclusion and further work}

In this paper, we introduce a game analogue of the generalized graph Tur\'an problem, the Constructor-Blocker game, and provide results for several natural choices for $H$ and $F$, but there are many more that remain unexplored. It would be interesting to see where the score of the game lies for a variety of other combinations both $H$ and $F$, and what is its relation to the corresponding generalized Tur\'an number. We are particularly curious about the leading term of $g(n,P_4,P_5)$, for which we provide bounds in Theorem~\ref{p4p5}.

As a natural extension of our game setup we could allow $H$ and/or $F$ to be representatives of a graph-theoretic structure that depends on the order of the board $n$, say a spanning graph like a spanning tree or a Hamilton cycle.

Also, what we study is sometimes referred to as an \emph{unbiased game}, where each player claims one edge per move. It would be interesting to study the $(1:b)$ biased version, for an integer $b>1$, in which Blocker claims $b$ edges per move.

In an effort to design a game version of the generalized graph Tur\'an problem, we opted for combining aspects of the Maker-Breaker positional games and the saturation games in the way described in the introduction. The same general approach readily gives several more possibilities for the game definition. In order to describe them we will refer to the edges claimed by the first player (resp.~second player) as red (resp.~blue) edges, and all the edges claimed by both players (i.e.~the union of red and blue edges) as black edges.
\begin{itemize}
    \item[(a)] Adhering more to the saturation games, we could forbid moves of both players that create a copy of $F$ in the black graph, counting the copies of $H$ also in the black graph at the end of the game.
    
    \item[(b)] In another approach somewhere in between (a) and our Constructor-Blocker games, we could forbid moves (of both players) that create a copy of $F$ in the black graph, but count the copies of $H$ in the red graph at the end of the game.
    
    \item[(c)] Then, there is an option of forbidding each player to claim a copy of $F$ in his own color, while counting the copies of $H$ in the black graph at the end of the game.
\end{itemize}
In all three versions the first player would try to maximize the number of counted copies of $H$, while the second player's goal would be to minimize it. As in the original Constructor-Blocker game, we are curious to know what the score of the game would be for various graph pairs in each of these cases, and where would it stand related to the corresponding generalized Tur\'an number, as well as the score in other versions.

\section*{Acknowledgments}
We would like to thank the anonymous referees whose comments and suggestions improved our paper.

\end{document}